\documentclass[12pt,a4paper]{article}
\usepackage{latexsym,amsfonts,amssymb,amsmath,longtable,amsthm,mathrsfs}
\usepackage[unicode,final,hyperindex]{hyperref}
\renewcommand{\le}{\leqslant}
\renewcommand{\ge}{\geqslant}

\newcommand{\No}{No}

\newtheorem{Lemma}{{\bfseries Lemma}}

\newtheorem{Theo}[Lemma]{{\bfseries Theorem}}
\theoremstyle{definition}

 \DeclareMathOperator{\Aut}{Aut}

\begin{document}

\begin{center}
 {\bfseries On the structure of groups, possessing Carter subgroups of odd
order}

E.P.Vdovin
\end{center}

\begin{flushright}

To Yuriy Leonidovich Ershov

on the occasion of his 75th birthday
 \end{flushright}

{\small Abstract. In the note we prove that all composition factors of a finite group possessing a Carter subgroup of odd order either are abelain, or are isomorphic to  $L_2(3^{2n+1})$.}


\section*{Introduction}

A known result by Glauberman and Thompson states, that a finite simple group
can not includes a self-normalizing Sylow  $p$-subgroup for
$p\ge5$ (see \cite[Theorem X.8.13]{HupBlackburn}, for example). Later, in
\cite[Corollary 1.2]{GuralnkickMalleNavarro} Guralnick, Malle, and Navarro
obtain a generalization of this result, proving that in any simple group  $G$
for a Sylow subgroup  $P$ of odd order the equality
$$N_G(P)=PC_G(P)$$ can not be fulfilled. This result is obtained by the authors
as a corollary to the following theorem.

\begin{Theo}\label{GurMalNavMain} {\em
\cite[Theorem~1.1]{GuralnkickMalleNavarro}}
Let $p$ be an odd prime and $P$ a Sylow $p$-subgroup of the finite
group  $G$. If  $p=3$, assume that   $G$ has no composition factors of type
$L_2(3^f)$, $f=3^a$ with $a\ge 1$.
\begin{itemize}
 \item[{\em (1)}] If $P=N_G(P)$, then $G$ is solvable.
 \item[{\em (2)}] If $N_G(P)=PC_G(P)$, then $G/O_{p'}(G)$ is solvable.
\end{itemize}
\end{Theo}

In the paper we prove a generalization of the first statement of the theorem.

\begin{Theo}\label{Main} {\em (Main Theorem)}
Assume that  $G$ possesses a Carter subgroup of odd order, Then each
composition factor of  $G$ either is abelian, or is isomorphic to
$L_2(3^{2n+1})$, $n\ge1$. Moreover, if  $3$ does not divide the order of a
Carter subgroup, then $G$ is solvable.
\end{Theo}

Clearly, item (1) of Theorem  \ref{GurMalNavMain} follows from Lemmas
\ref{CarterInInducedGroup} and \ref{CarterOfOddOrderAlmostSimple} (see the
proof in the end of the paper).

\section{Notations}

In the paper only finite groups are considered, so the term ``group'' is
always used in the meaning ``finite group''.

The notation in the paper agrees with that of~\cite{ATLAS}. Recall that a nilpotent selfnormalizing subgroup is called a
{\em Carter subgroup}. A non-refinable normal series of a group is called a {\em chief series}. A composition series is called an {\em $(rc)$-series}\footnote{this term is introduced by V.A.Vedernikov in \cite{Vedernikov1}}, if it is a refinement of a chief series.

Let $A,B,H$ be subgroups of  $G$ such that $B$  is normal in $A$. Define  $N_H(A/B):=N_H(A)\cap
N_H(B)$ to be  the {\em normalizer} of $A/B$ in $H$. If $x\in N_H(A/B)$, then  $x$  induces an automorphism on
$A/B$ acting by $Ba\mapsto B x^{-1}ax$. Thus there exists a homomorphism  $N_H(A/B)\rightarrow \Aut(A/B)$. The image of 
$N_H(A/B)$ under the homomorphism is denoted by  $\Aut_H(A/B)$ and is called the {\em group of $H$-induced automorphisms} of  $A/B$, while the kernel of the homomorphism is denoted by  $C_H(A/B)$ and is called the {\em centralizer of $A/B$ in $H$}. If  $B=1$, then we use the notation $\Aut_H(A)$. Notice that   $\Aut_G(A)$ sometimes  is called the  automizer of  $A$ in $G$. Groups of induced automorphisms are introduced by F.Gross in  \cite{GrossExistence}, where the author says that this notion is taken from unpublished Wielandt's lectures. Evidently,  $C_H(A/B)=C_G(A/B)\cap H$, so
$$\Aut_H(A/B)=N_H(A/B)/C_H(A/B)\simeq N_H(A/B)C_G(A/B)/C_G(A/B)\le
\Aut_G(A/B),$$ i.e.  $\Aut_H(A/B)$ can be naturally considered as a subgroup of $\Aut_G(A/B)$, and we think of  $\Aut_H(A/B)$ as a subgroup of 
$\Aut_G(A/B)$ without additional clarifications.

We need the following result.

\begin{Lemma}\label{GeneralizedJordanHolder} {\em \cite[Theorem 1]{VdoMalt}
(Generalized Jordan-H\"{o}lder theorem)}
Let  $$G=G_0\supset G_1\supset\ldots\supset G_n=1$$ 
be an $(rc)$-series of $G$, denote $G_{i-1}/G_i$ by  $S_i$.
Assume that  $$G=H_0\supset
H_1\supset\ldots\supset H_n=1$$
is a composition series of $G$ and $T_i=H_{i-1}/H_i$. Then there exists a permutation  $\sigma\in Sym_n$ such that for every section $T_i$ the inclusion
$\Aut_G(T_i)\le\Aut_G(S_{i\sigma})$ holds. Moreover, if the second series is also an $(rc)$-series, then   $\sigma$ can be chosen so that the isomorphims $\Aut_G(T_i)\simeq\Aut_G(S_{i\sigma})$ holds.
\end{Lemma}

\section{Proof of the main theorem}

We divide the proof of the main theorem into several lemmas.

\begin{Lemma}\label{CarterInInducedGroup}
Let  $K$ be a Carter subgroup of  $G$ and $$G=G_0\supset
G_1\supset\ldots\supset
G_n=1$$  be an
$(rc)$-series of  $G$. Then for every nonabelian composition factor 
$S$ of $G$ there exists  $i$ such that  $G_{i-1}/G_i\simeq S$ and 
$\Aut_K(G_{i-1}/G_i)$ is a Carter subgroup of
$\Aut_G(G_{i-1}/G_i)$.
\end{Lemma}

\begin{proof}
The claim follows by induction on the length of the chief series, whose refinement is the  $(rc)$-series, and  \cite[Lemma~3]{CartExist}.
\end{proof}

\begin{Lemma}\label{CarterOfOddOrderAlmostSimple} {\em (mod CFSG)}
Let  $G$ be a finite almost simple group, possessing a Carter subgroup  
$K$ of odd order. Then  $G\simeq L_2(3^{2n+1})\leftthreetimes\langle
\varphi\rangle$, where $n\ge1$ and $\varphi$ is a field automorphism of $G$
of order~$2n+1$.

In particular, if a Sylow   $3$-subgroup of  $G$  is a Carter subgroup, then $G\simeq L_2(3^{3^n})\leftthreetimes\langle
\varphi\rangle$, where $n\ge1$ and $\varphi$ is a field automorphism of   $G$ of order~$3^n$.
\end{Lemma}

\begin{proof}
The clain follows from the classification of Carter subgroups given in \cite[Tables 7--10]{CartMathTr}. Notice that only this lemma in the paper uses the classification of finite simple groups.
\end{proof}

Now we are ready to proof the main result of the paper (Theorem \ref{Main}).
Assume that a finite group  $G$ possesses a Carter subgroup  $K$ of odd order. Assume that there exists a nonabelian composition factor $S$ of $G$. Then by Lemma  \ref{CarterInInducedGroup}, there exist subgroups 
$A,B$ of $G$ such that   $A/B\simeq S$ and $\Aut_K(A/B)$ is a Carter subgroup of  $\Aut_G(A/B)$. By Lemma \ref{CarterOfOddOrderAlmostSimple}
we obtain $S\simeq L_2(3^{2n+1})$. Notice that by \cite[Table~10]{CartMathTr} it follows that in this case $|\Aut_K(A/B)|$ is divisible by $3$, i.e. $|K|$ is divisible by $3$ as well. Therefore, if  $|K|$ is not divisible by  $3$, then  $G$ is solvable.

Notice that statement  (1) in Theorem \ref{GurMalNavMain} can be obtained by exactly the same arguments.


\begin{thebibliography}{99}
\bibitem{HupBlackburn} {\itshape B. Huppert and N. Blackburn}, Finite Groups
III, Springer-Verlag, Berlin, New York, 1982.

\bibitem{GuralnkickMalleNavarro} {\itshape R.M. Guralnick, G. Malle,
G. Navarro}, Self-normalizing Sylow subgroups,  Proc. Amer. Math. Soc.,
{\bfseries 132}, \No. 4 (2004), 973–979.

\bibitem{ATLAS}
{\itshape J.H. Conway, R.T. Curtis, S.P. Norton, R.A. Parker, R.A. Wilson},
Atlas of Finite Groups, Clarendon Press, Oxford, 1985.

\bibitem{Vedernikov1} {\itshape V.A.Vedernikov}, Finite groups with Hall 
$\pi$-subgroups, Math.Notes, {\bfseries 203} (2012), \No 3, 326--350.

\bibitem{GrossExistence}
{\itshape F. Gross}, On the existence of Hall subgroups, J. Algebra,
{\bfseries 98}, \No 1 (1986), 1--13.

\bibitem{VdoMalt} {\itshape E.P.Vdovin}, Groups of induced automorphisms
and their application to studying the existence problem for Hall subgroups,
Algebra and logic, {\bfseries 53} (2014), \No 5, 418--421.

\bibitem{CartExist} {\itshape E.P.Vdovin}, On the conjugacy problem for Carter subgroups, Siberian
Mathematical Journal, {\bfseries 47} (2006), \No~4, 597--600.

\bibitem{CartMathTr} {\itshape E.P.Vdovin},
Carter subgroups of finite groups, Siberian Advances in Mathematics, {\bfseries 19}
(2009), \No~1, 24--74.
(corrigendum:
\verb{http://math.nsc.ru/~vdovin/Papers/carter_eng_corrigendum.pdf{)

\end{thebibliography}
\end{document}